\documentclass[12pt, reqno]{amsart}

\usepackage{amssymb}
\usepackage{amsbsy}
\usepackage{leftidx}

\newcommand{\orbfont}[1]{\mathfrak{#1}}
\newcommand{\grfont}[1]{\mathcal{#1}}

\newcommand{\homeo}{\cong}
\newcommand{\isom}{\cong}

\newcommand{\homotopic}{\simeq}
\newcommand{\weakhtpic}{\sim}
\newcommand{\equivalent}{\sim}
\newcommand{\disunion}{\coprod}
\newcommand{\smashprod}{\wedge}

\newcommand{\piorb}{\pi^{\mathrm{orb}}}
\newcommand{\piorbext}{\hat{\pi}^{\mathrm{orb}}}
\newcommand{\sthorb}{\stho^{\mathrm{orb}}}
\newcommand{\directsum}{\bigoplus}
\newcommand{\Hom}{\operatorname{Hom}}

\newcommand{\im}{\operatorname{im}}
\newcommand{\map}{\operatorname{map}}
\newcommand{\fibover}[2]{\mathbin{\leftidx{_{#1}}{\times}{_{#2}}}}

\newcommand{\catfont}[1]{\mathtt{#1}}
\newcommand{\Orb}{\catfont{Orb}} 
\newcommand{\Gpd}{\catfont{Gpd}} 

\newcommand{\bbR}{\mathbb{R}}

\newcommand{\bbS}{\mathbb{S}}

\newcommand{\stho}{\varpi}
\newcommand{\semidirect}{\ltimes}

{\theoremstyle{plain}
\newtheorem{theorem}{Theorem}[section]
\newtheorem{lemma}[theorem]{Lemma}
\newtheorem{corollary}[theorem]{Corollary}

\newtheorem{proposition}[theorem]{Proposition}
}

{\theoremstyle{definition} 
\newtheorem{remark}[theorem]{Remark}
\newtheorem{exa}[theorem]{Example}

}

\newenvironment{example}{\begin{exa}}{\hfill $\Diamond$\end{exa}}

\usepackage{mathrsfs}
\usepackage[all]{xy}
\usepackage{amssymb}
\usepackage{graphicx}
\usepackage[margin=1.5in]{geometry}

\newtheorem*{maintheorem}{Theorem \ref{thm::maintheorem}}
\newtheorem*{mainprop}{Proposition \ref{prop::mainprop}}
\newtheorem*{maincor}{Corollary \ref{cor::maincor}}
\newcommand{\Gbar}{\underline{G}}

\title{Orbifolds and Stable Homotopy Groups}
\author{Johann K. Leida}
\address{Dept. of Mathematics\\
University of Wisconsin\\
Madison, WI\\
USA}
\email{leida@math.wisc.edu}
\date{\today}

\thanks{The author was supported by an NSF Graduate Research
Fellowship.}

\begin{document}

\begin{abstract}
Lie groupoids generalize transformation groups, and so provide a
natural language for studying orbifolds
\cite{moerdijk;orbifolds-groupoids-introduction} and other
noncommutative geometries. In this paper, we investigate a connection
between orbifolds and equivariant stable homotopy theory using such
groupoids. A different sort of twisted sector, along with a classical
theorem of tom Dieck \cite{tomdieck;transformation-groups;;1987},
allows for a natural definition of \emph{stable orbifold homotopy
groups}, and motivates defining \emph{extended unstable orbifold
homotopy groups} generalizing previous definitions.
\end{abstract}

\maketitle

\section{Introduction}\label{sec::intro}

Homotopy groups for orbifolds have been defined in a variety of ways,
going back to Thurston's orbifold fundamental group
\cite{thurston-notes} and culminating in Moerdijk's elegant treatment,
which we briefly recall here. We assume the reader is familiar with
the category $\Gpd$ of orbifold groupoids and their homomorphisms, as
well as Morita equivalences, classifying spaces, groupoid actions, and
other basic notions--see
\cite{moerdijk;orbifolds-groupoids-introduction} or \cite{orbibook}
for an introduction and further references. By an \emph{orbifold
groupoid} $\grfont{G}$, we mean a proper foliation Lie groupoid. We
write $G_{0}$ for the manifold of objects and $G_{1}$ for the manifold
of arrows. Let $\grfont{G}$ be an orbifold groupoid, and let $x\in
G_{0}$ be a base point. The \emph{$n^{\text{th}}$ orbifold homotopy
group} of $\grfont{G}$ based at $x$ is
\begin{equation}\label{}
\piorb_{n}(\grfont{G}, x):=\pi_{n}(B\grfont{G}, [x]),
\end{equation}
where $B\grfont{G}$ is the classifying space of $\grfont{G}$, and
$[x]$ is the point in $B\grfont{G}$ corresponding to the object $x$.

These groups have several agreeable properties. They are \emph{Morita
invariant}, so that they descend to the localized category $\Orb$ of
orbifolds. They also generalize Thurston's fundamental group, and
agree with alternative definitions, such as
\cite{chen;homotopy-theory-orbispaces}, in higher degrees. On the
other hand, equivariant homotopy theory reveals that these groups are
insufficient, for if $\grfont{G}=G\semidirect M$ is a translation
groupoid corresponding to the quotient orbifold $\orbfont{X}=M/G$,
then $B\grfont{G}$ is homotopic to the Borel construction
$EG\times_{G}M$ (see Appendix \ref{app::borelconstruction}), which
fails to capture the $G$-homotopy type of $M$.
\begin{example}\label{ex::Disk}
Let $D$ be a disk with a smooth, fixed-point free action of the
icosahedral group $\mathcal{I}$. Such an action is described in
\cite[pp. 55--58]{bredon;compact-transformation-groups}, and was first
constructed simplicially by Floyd and Richardson
\cite{floyd-richardson;action-finite-group-without-stationary}. Then
the map $f:D\rightarrow \{\text{pt}\}$ is an equivariant map that is a
nonequivariant homotopy equivalence. Accordingly, the induced map
$f:E\mathcal{I}\times_{\mathcal{I}}D\rightarrow B\mathcal{I}$ is a
homotopy equivalence. However, $\mathcal{I}\semidirect D$ is certainly
not the same orbifold as $\mathcal{I}\semidirect \{\text{pt}\}$, as
the former has no point with isotropy $\mathcal{I}$.
\end{example}
For a less complicated example, take any group $G$. Then the unique
$G$-map $EG\rightarrow \{\text{pt} \}$ also induces a homotopy
equivalence between the Borel constructions. However, here we must
drop some smoothness or properness conditions; in other words, at
least one of these translation groupoids is not a (finite dimensional)
orbifold. In any case, one ought not call the translation groupoids in
either example ``homotopy equivalent.''

These examples suggest that it is necessary to locate additional fixed
point data in topological groupoids. For in the equivariant situation,
one would consider the homotopy types of various fixed point sets and
easily distinguish the spaces above. Alternatively, (integer graded)
stable equivariant homotopy groups also do the trick, given tom
Dieck's isomorphism \cite{tomdieck;transformation-groups;;1987}:
\begin{equation}\label{eqn::tomDieck}
\stho_{n}^{G}(X)\isom
\directsum_{(H)}\stho^{W_{G}H}_{n}(EW_{G}H_{+}\smashprod X^{H}).
\end{equation}
Here $G$ is supposed to be a compact Lie group acting on the pointed
$G$-space $X$, the sum runs over all conjugacy classes of subgroups,
$W_{G}H=N_{G}H/H$, and $+$ denotes a disjoint $G$-fixed base
point. This result demonstrates that these groups depend mainly on
fixed point data and isotropy groups, so that one might hope that they
are ``Morita invariant'' in some sense. We make this hope precise by
defining intrinsic \emph{stable orbifold homotopy groups}
\begin{equation*}
\sthorb_{n}(\orbfont{X}) := \stho_{n}(B\widetilde{\grfont{G}}_{+}),
\end{equation*}
where $\grfont{G}$ is an orbifold groupoid corresponding to
$\orbfont{X}$ and $\widetilde{\grfont{G}}$ is the groupoid of \emph{fixed
point sectors} defined in Section \ref{sec::sectors}. We can also
define \emph{extended unstable orbifold homotopy groups} 
\begin{equation*}
\piorbext_{n}(\orbfont{X}, (x,H)) := \pi_{n}(B\widetilde{\grfont{G}}, [x,H])
\end{equation*}
where $(x,H)$ is a base point in the fixed point sectors. The main
theorem is then:
\begin{maintheorem}
The stable orbifold homotopy groups $\sthorb_{n}(\orbfont{X})$ and
extended orbifold homotopy groups $\piorbext_{n}(\orbfont{X}, (x, H))$
are orbifold invariants. 
\end{maintheorem}

We apply this result to prove the following proposition in the case
where $\orbfont{X}$ is a quotient.
\begin{mainprop}
Let $\orbfont{X}=M/G$ where $M$ is a smooth manifold and $G$ is a
compact Lie group acting smoothly and almost freely. Then the total
stable equivariant homotopy group 
\begin{equation*}
\stho^{G}_{\mathrm{tot}}(M_{+}):=\bigoplus_{n} \stho^{G}_{n}(M_{+})
\end{equation*}
is an orbifold invariant.
\end{mainprop}
The proof amounts to a calculation identifying the stable equivariant
homotopy group with the orbifold stable homotopy group in this special
case. This is analogous to the situation with orbifold $K$-theory
\cite{adem-ruan;twisted-orbifold-K-theory}; there, one sees that
$K_{G}^{*}(M)\isom K_{G'}^{*}(M')$ whenever $M/G$ and $M'/G'$ are the
same orbifold by identifying each with the intrinsically defined
orbifold $K$-theory: 
\begin{equation*}
\xymatrix{ 
K_{G}^{*}(M) & K_{\text{orb}}^{*}(\orbfont{X})\ar[l]_{\isom}
\ar[r]^{\isom} & K_{G'}^{*}(M').  
}
\end{equation*}
In fact, in the very special case of finite group quotients (so-called
\emph{global quotients}), we have:
\begin{maincor}
If $M/G$ and $M'/G'$ are two global quotient presentations of the same
orbifold $\orbfont{X}$, then there are isomorphisms
\begin{equation*}
\xymatrix{ 
\stho^{G}_{n}(M_{+}) & \sthorb_{n}(\orbfont{X})
\ar[l]_{\isom}\ar[r]^{\isom} & \stho_{n}^{G'}(M_{+}') 
}
\end{equation*}
for each integer $n$.
\end{maincor}
We emphasize that the stable homotopy calculation requires information
about fixed points for all subgroups, whereas the $K$-theory is
detected on fixed point sets for group elements (or, equivalently, on
cyclic subgroups).

We make our definitions within the smooth category, since that is
where our present applications lie. However, the reader will note that
most of the definitions go through for well-behaved topological
groupoids. The first order of business is to organize the fixed point
data of an orbifold using its groupoid presentations. This is
accomplished in \S\ref{sec::sectors}. Fixed point data in hand, we
define the novel unstable and stable homotopy groups in
\S\ref{sec::homotopyGroups}. After seeing how the new invariants
generalize classical ones, we close with some ideas for future
directions and applications in \S\ref{sec::future}.

The author would like to thank Peter May and Christopher French for
several helpful conversations, W.~G.~Dwyer for suggesting the
Floyd-Richardson example, and also Ieke Moerdijk for many helpful
comments on an earlier version of this paper.

\section{Sectors}\label{sec::sectors}

For an orbifold groupoid $\grfont{G}$, we need to understand the fixed
point data $\grfont{G}$ encodes, and then show that it is Morita
invariant. Recall that the \emph{isotropy group} $G_{x}$ of an object
$x\in G_{0}$ is the group $s^{-1}(x)\cap t^{-1}(x)$ of all arrows that
start and end at $x$. Let
\begin{equation}\label{}
\widetilde{S}(\grfont{G}):=\{(x,H) \mid x\in G_{0}, H\subseteq G_{x} \}
\end{equation}
be the set of all subgroups of the groupoid $\grfont{G}$--that is, the
set of all the subgroups of all the isotropy groups. The goal is to
correctly topologize this set of ``fixed points.'' The correct
topology should specialize to a disjoint union of ordinary fixed point
sets in the case of a translation groupoid $G\semidirect M$. In
particular, one expects that for each group $H$, the subset
$\widetilde{S}^{H}(\grfont{G}):=\{(x,K) \mid K\isom H \}$ will be open and
closed, so that $\widetilde{S}^{H}(\grfont{G})$ is a union of connected
components. It then suffices to topologize each
$\widetilde{S}^{H}(\grfont{G})$ separately. These \emph{fixed point
sectors} are closely related to the twisted multisectors appearing in
Chen-Ruan cohomology and orbifold $K$-theory. In fact, both
constructions are rooted in Kawasaki's earlier work
\cite{kawasaki-signature-theorem-V-manifolds}.

\subsection{Definitions}\label{subsec::sectorDefns}

Assume that $\grfont{G}$ is an orbifold groupoid; in particular, that
$\grfont{G}$ has all isotropy groups $G_{x}$ finite. Then for each
finite group $H$ of order $k$, we can identify
$\widetilde{S}^{H}(\grfont{G})$ with a subquotient space of the
\emph{$k$-sectors} $\widetilde{S}^{k}(\grfont{G})$. Recall that elements
of $\widetilde{S}^{k}(\grfont{G})$ are ordered tuples
$(g_{1},\dots,g_{k})$ of arrows in $\grfont{G}$ that all begin and end
at the same place. More formally, the $k$-sectors are topologized as
the fibered product
\begin{equation*}
\xymatrix{
\widetilde{S}^{k}(\grfont{G})\ar[r]\ar[d] & (G_{1})^{k}\ar[d]^{(s,t)^k}\\
G_{0}\ar[r]^-{\Delta^{k}}      & (G_{0}\times G_{0})^{k}
}
\end{equation*}
where $\Delta$ is the diagonal map. Note that
$\widetilde{S}^{k}(\grfont{G})$ is a smooth manifold when $\grfont{G}$ is
a Lie groupoid.

Now let $\widetilde{P}^{H}(\grfont{G})$ be the subset of those
$(g_{1},\dots,g_{k})\in \widetilde{S}^{k}(\grfont{G})$ such that the
$g_{i}$ are all distinct and, as a subset of $G_{1}$, form a group
isomorphic to $H$. The symmetric group $\mathfrak{S}_{k}$ acts freely
on $\widetilde{P}^{H}(\grfont{G})$ by permuting coordinates, and the
quotient is in bijection with $\widetilde{S}^{H}(\grfont{G})$. Topologize
$\widetilde{S}^{H}(\grfont{G})$ by declaring this to be a
homeomorphism. Finally, let
\begin{equation}
\widetilde{S}(\grfont{G}):=\disunion_{H} \widetilde{S}^{H}(\grfont{G})
\end{equation}
as $H$ runs over representatives of all isomorphism classes of
subgroups of $\grfont{G}$. In many cases of interest (including
compact orbifolds and global quotient orbifolds), $\grfont{G}$ will
only have finitely many isomorphism types of subgroups, so that the
union is taken over a finite set.

\begin{remark}
In fact, $\tilde{S}(\grfont{G})$ has a natural smooth structure. We
will see this later in the proper \'etale case; the extension to the
proper foliation case is not difficult.
\end{remark}

So, is this topology ``correct?'' Consideration of a few examples
leads us to believe that it is.
\begin{example}
Let $\grfont{G}=M$ be a manifold viewed as a trivial or \emph{unit}
orbifold groupoid. Then
$\widetilde{S}(\grfont{G})=\widetilde{S}^{\langle
1\rangle}(\grfont{G})$ consists of only the \emph{trivial sector}
corresponding to the trivial group $\langle 1\rangle$, and
$\widetilde{S}^{\langle 1\rangle}(\grfont{G})=\widetilde{P}^{\langle
1\rangle}(\grfont{G})/\mathfrak{S}_{1}\homeo M$.
\end{example}
\begin{example}
Let $\grfont{G}=G$ be a finite group. Then $\widetilde{S}(\grfont{G})$ is the set of
all subgroups of $G$ endowed with the discrete topology. 
\end{example}
\begin{example}\label{ex::transSectors}
Both previous examples are special cases of translation groupoids. Let
$G$ be a Lie group acting smoothly and almost freely on the manifold
$U$, so that $\grfont{G}=G \semidirect U$ is an orbifold
groupoid. Then $\widetilde{S}(\grfont{G}) = \disunion_{H\subseteq
G}U^{H}$ is the disjoint union of the fixed point sets. To see this,
recall that in this situation
\begin{equation*}
\widetilde{S}^{k}(\grfont{G}) = \disunion_{(g_{1},\dots, g_{k})\subseteq
G^{k}} U^{g_{1}}\cap\dots\cap U^{g_{k}}.
\end{equation*}
So we have 
\begin{equation*}
\widetilde{P}^{H}(\grfont{G}) = \disunion_{\substack{(g_{1},\dots,
g_{k})=K\subseteq G\\
K\isom H}} U^{K},
\end{equation*}
where $U^{K}$ is appearing $|K|!$ times, and hence 
\begin{equation}
\widetilde{S}^{H}(\grfont{G})=\widetilde{P}^{H}(\grfont{G})/\mathfrak{S}_{k} =
\disunion_{\substack{K\subseteq G\\
K\isom H}} U^{K}.
\end{equation}
\end{example}

Confident in this topology for $\widetilde{S}(\grfont{G})$, we want to add
some arrows and turn it into a groupoid. After all, the fixed point
sets $U^{K}$ in the last example come equipped with natural actions of
the normalizers $N_{G}K$. Similar information may be extracted
from more general groupoids. In fact, this is easy; the well-known
conjugation action of $\grfont{G}$ on the $k$-sectors induces an
action on $\widetilde{S}^{H}(\grfont{G})$.
\begin{lemma}
Let $\grfont{G}$ be an orbifold groupoid. Then for each $H$, the fixed
point sector $\widetilde{S}^{H}(\grfont{G})$ is naturally a smooth
$\grfont{G}$-space. Thus, the associated translation groupoid,
$\widetilde{\grfont{G}}^{H}:=\grfont{G}\semidirect
\widetilde{S}^{H}(\grfont{G})$, is again an orbifold groupoid.
\end{lemma}
\begin{proof}
The action of $\grfont{G}$ on $\widetilde{S}^{k}(\grfont{G})$ is given
as follows: the anchor map $\pi :\widetilde{S}^{k}(\grfont{G})\to
G_{0}$ sends $(g_{1},\dots,g_{k})$ to their common source/target. The
action map $\mu :G_{1}\fibover{s}{\pi}\widetilde{S}^{k}(\grfont{G})\to
\widetilde{S}^{k}(\grfont{G})$ is given by conjugation:
$h(g_{1},\dots, g_{k})=(hg_{1}h^{-1},\dots ,hg_{k}h^{-1})$. These maps
are smooth, making the $k$-sectors into a smooth
$\grfont{G}$-space. So all that remains is to observe that
$\pi|_{\widetilde{P}^{H}(\grfont{G})}$ is $\mathfrak{S}_{k}$-invariant
and that the restriction of $\mu$ is equivariant with image in
$\widetilde{P}^{H}(\grfont{G})$.
\end{proof}

We write $\widetilde{\grfont{G}}$ for $\grfont{G}\semidirect
\widetilde{S}(\grfont{G})$, and call it the groupoid of \emph{fixed point
sectors}. Its subgroupoid $\widetilde{\grfont{G}}^{H}$ is called the
\emph{$H$-fixed sector}. We can now complete our earlier translation
groupoid example.
\begin{example}
Let $\grfont{G}=G\semidirect U$ as before. Then the identification
of Example \ref{ex::transSectors} extends to a groupoid isomorphism
\begin{equation}\label{}
\widetilde{\grfont{G}}\isom G \semidirect (\disunion_{K\subseteq G} U^{K}).
\end{equation}
Further, one readily sees that the inclusion 
\begin{equation*}
\disunion_{(K)\subseteq G} N_{G}K\semidirect U^{K}\hookrightarrow G \semidirect (\disunion_{K\subseteq G} U^{K})
\end{equation*}
is a weak equivalence\footnote{Moerdijk calls such homomorphisms
\emph{equivalences},
c.f. \cite[p. 209]{moerdijk;orbifolds-groupoids-introduction}. This
could be misleading, since they do not form a symmetric relation.},
where the union on the left now runs over conjugacy classes of
subgroups. So, in particular, $\widetilde{G \semidirect U}$ is Morita
equivalent to $\disunion_{(K)\subseteq G} N_{G}K\semidirect
U^{K}$. This example will be quite useful when we study the fixed
point sectors of proper \'etale groupoids.
\end{example}
Note that the quotient space $|\widetilde{\grfont{G}}|$ is the set of
points 
\begin{equation}\label{}
\{(x,(H)_{G_{x}})\mid x\in |\grfont{G}|, H\subseteq G_{x}\},
\end{equation}
where $(H)_{G_{x}}$ indicates the conjugacy class of $H$ in
$G_{x}$. This recovers Kawasaki's description.

In considering the fixed point set $M^{H}$ of a group action, one
often disregards the trivial $H$ action and instead focuses on the
action of the Weyl group $WH=NH/H$. We can also do this in the
groupoid case. Let $\grfont{G}$ be an orbifold groupoid. We consider
the following subset $\grfont{K}_{H}$ of arrows in
$\widetilde{\grfont{G}}^{H}$:
\begin{equation}
\grfont{K}_{H}:=\{(l, (x,L)) \mid l\in L\}.
\end{equation}
Then $\grfont{K}_{H}$ forms a wide normal subgroupoid of
$\widetilde{\grfont{G}}^{H}$. Recall that \emph{wide} simply means that a
subgroupoid contains all identity arrows, and \emph{normal} means
that conjugates of arrows in the subgroupoid are again in the
subgroupoid, whenever conjugation makes sense. $\grfont{K}_{H}$ is
also \emph{totally intransitive}, in that it consists entirely of
isotropy arrows.

Define $\overline{\grfont{G}}^{H} := \widetilde{\grfont{G}}^{H} /
\grfont{K}_{H}$, and let $\overline{\grfont{G}} :=
\disunion_{H}\overline{\grfont{G}}^{H}$, where $H$ runs over
isomorphism classes of subgroups as before. These quotient groupoids
have the same set of objects as the fixed point sectors, but fewer
arrows. In fact, one can check that for $\grfont{G}=G\semidirect M$,
we have $\overline{G\semidirect M}$ Morita equivalent to
$\disunion_{(H)} W_{G}H\semidirect M^{H}$. In general, the quotient
space $|\overline{\grfont{G}}|$ is homeomorphic to
$|\widetilde{\grfont{G}}|$, but has a different orbifold structure
(with less isotropy). $\overline{\grfont{G}}$ is called the groupoid
of \emph{reduced fixed point sectors}.

The fixed point sectors $\widetilde{\grfont{G}}$ (and the reduced
sectors $\overline{\grfont{G}}$) are natural in $\grfont{G}$, inherit
many of its properties, and respect Morita equivalence. A Lie groupoid
homomorphism $\phi:\grfont{G}\rightarrow \grfont{H}$ is
\emph{faithful} if it is faithful as a functor (i.e., injective on
$\Hom$ sets).
\begin{lemma}\label{lem::sectors}
The fixed point sector construction $\widetilde{\grfont{G}}$ and its
reduced version $\overline{\grfont{G}}$ are both functorial with
respect to faithful orbifold groupoid homomorphisms, and both respect
Morita equivalences. Moreover, $\widetilde{\grfont{G}}$ and
$\overline{\grfont{G}}$ are orbifold groupoids, and are \'etale if and
only if $\grfont{G}$ is such a groupoid.
\end{lemma}

\begin{proof}
Suppose $\phi :\grfont{G}\rightarrow \grfont{H}$ is a homomorphism
between two orbifold groupoids. Then $\phi$ induces a homomorphism
$\phi_{*}:\widetilde{S}^{k}(\grfont{G})\to\widetilde{S}^{k}(\grfont{H})$ for
each $k$, where
\begin{equation*}
\phi_{*}(g_{1},\dots,g_{k}):=(\phi (g_{1}),\dots,\phi(g_{k})).
\end{equation*}
If $\phi$ is faithful, then for each $H$ the restriction of $\phi_{*}$
to $\widetilde{P}^{H}(\grfont{G})$ lands in
$\widetilde{P}^{H}(\grfont{H})$, and is
$\mathfrak{S}_{k}$-equivariant. Thus there are induced maps
$\widetilde{\phi}^{H}:\widetilde{S}^{H}(\grfont{G}) \to
\widetilde{S}^{H}(\grfont{H})$, with union
$\widetilde{\phi}:\widetilde{S}(\grfont{G}) \to
\widetilde{S}(\grfont{H})$. To extend $\widetilde{\phi}$ to arrows, we
just restrict the map $\phi_{1}\times \widetilde{\phi}_{0}:G_{1}\times
\widetilde{S}(\grfont{G})\to H_{1}\times \widetilde{S}(\grfont{H})$ to
$G_{1}\fibover{s}{\pi}\widetilde{S}(\grfont{G})$. The image is
automatically in $H_{1}\fibover{s}{\pi}\widetilde{S}(\grfont{H})$, and
it is easy to check that this is a groupoid homomorphism. The kernel
$\grfont{K}_{H}$ is preserved under this homomorphism, so we also
obtain an induced homomorphism between the reduced sectors.

Next suppose $\phi$ is a weak equivalence. We show $\widetilde{\phi}$
is as well, so that $\widetilde{\grfont{G}}$ and
$\widetilde{\grfont{H}}$ are Morita equivalent whenever $\grfont{G}$
and $\grfont{H}$ are. This is actually a special case of a more
general situation. If $E$ is a $\grfont{H}$-space, then the pullback
$\phi^{*}E:=E\fibover{\pi}{\phi}G_{0}$ is a $\grfont{G}$-space, and
$\phi$ induces a homomorphism between the associated translation
groupoids $\grfont{G}\semidirect \phi^{*}E$ and $\grfont{H}\semidirect
E$. When $\phi$ is a weak equivalence, this induced homomorphism is
automatically a weak equivalence. So it suffices to show that
$\widetilde{S}(-)$ is \emph{natural}, that is,
$\phi^{*}\widetilde{S}(\grfont{H})\isom \widetilde{S}(\grfont{G})$ as
$\grfont{G}$-spaces when $\phi$ is a weak equivalence.

It is enough to check this sector by sector, so we consider the map 
\begin{align*}
(\widetilde{\phi}_{H}, \pi ):\widetilde{S}^{H}(\grfont{G}) &\to \widetilde{S}^{H}(\grfont{H})\fibover{\pi}{\phi}G_{0}\\
[g_{1},\dots,g_{k}]&\mapsto ([\phi (g_{1}),\dots,\phi
(g_{k})],\pi[g_{1},\dots,g_{k}]).
\end{align*}
This is a smooth $\grfont{G}$-equivariant bijection, with a
smooth $\grfont{G}$-equivariant inverse given by
\begin{equation*}
([h_{1},\dots, h_{k}],x)\mapsto [\phi_{x}^{-1}(h_{1}),\dots, \phi_{x}^{-1}(h_{k})],
\end{equation*}
where $\phi_{x}^{-1}(h_{i})$ means the unique preimage in
$G_{x}$. Here we used the fact that $\phi$ is fully faithful to see
uniqueness, and the Cartesian square property in the definition of
weak equivalences
\cite[p. 209]{moerdijk;orbifolds-groupoids-introduction} to see that
the inverse is smooth. The proof for the reduced sectors is analogous.

Finally, we address the adjectives. First, $\grfont{G}$ is embedded in
$\widetilde{\grfont{G}}$ (and in $\overline{\grfont{G}}$) as the trivial
sector, so if $\widetilde{\grfont{G}}$ has any of the properties
below, they are automatically inherited by $\grfont{G}$. Conversely:
\begin{itemize}
\item If $\grfont{G}$ is proper, then so is $\widetilde{\grfont{G}}$. It's
enough to show the source map of $\widetilde{\grfont{G}}$ is proper, and
that is the map $\text{pr}_{2}:
G_{1}\fibover{s}{\pi}\widetilde{S}(\grfont{G})\rightarrow
\widetilde{S}(\grfont{G})$. This is proper because the source map $s$ of
$\grfont{G}$ is proper, and properness is stable under base change.
\item Foliation: the natural homomorphism
$\widetilde{\grfont{G}}\rightarrow \grfont{G}$ is easily seen to be
faithful. Hence all isotropy groups of $\widetilde{\grfont{G}}$ are
discrete when those of $\grfont{G}$ are.
\item The \'etale case is just like properness. Any translation
groupoid $\grfont{G}\semidirect E$ for an \'etale groupoid
$\grfont{G}$ is \'etale, since \'etale maps are preserved under base
change. 
\end{itemize}
The proof of these properties for $\overline{\grfont{G}}$ is
straightforward.
\end{proof} 

\begin{remark}
In fact, one expects these fixed point constructions to be functorial
with respect to all groupoid homomorphisms. Up to homotopy, this
certainly appears to be the case: we can replace
$\widetilde{\grfont{G}}$ with a topological category
$\catfont{Gps}/\grfont{G}$ for which functorality is clear. Here,
$\catfont{Gps}$ is the category of finite groups and surjective
homomorphisms, and $\catfont{Gps}/\grfont{G}$ has objects the groupoid
homomorphisms $H\stackrel{\rightarrow}{a} \grfont{G}$ for $H\in
\catfont{Gps}$, with arrows the commutative diagrams
\begin{equation*}
\xymatrix{
H \ar@{->>}[r]^{f} \ar[d]^{a} & K\ar[d]^{b}\\
\im(a)\subseteq \grfont{G} \ar[r]^{c_{g}} & \im(b)\subseteq\grfont{G},
}
\end{equation*}
where $c_{g}$ is conjugation by $g\in G_{1}$. The natural projection
$\catfont{Gps}/\grfont{G}\rightarrow \widetilde{\grfont{G}}$ induces a
homotopy equivalence of classifying spaces.
\end{remark}

\subsection{The Proper \'Etale Case}\label{subsec::etaleCase}    

Proper \'etale groupoids form a pleasant class of examples, and since
every orbifold $\orbfont{X}$ may be presented as such a groupoid, they
merit some elaboration. The key feature of proper \'etale groupoids is
that they are locally isomorphic to translation groupoids
$G_{x}\semidirect U_{x}$ for each $x\in G_{0}$. The fixed point sector
construction commutes with restriction, so in the proper \'etale case
the fixed point sectors $\widetilde{\grfont{G}}$ are themselves
locally isomorphic to $\widetilde{G_{x}\semidirect U_{x}}$. We have
seen in Example \ref{ex::transSectors} that this is Morita equivalent
to
\begin{equation*}
\disunion_{(H)}N_{G_{x}} H\semidirect U_{x}^{H}.
\end{equation*}
Similarly, $\overline{\grfont{G}}$ is locally modelled on 
\begin{equation*}
\disunion_{(H)}W_{G_{x}} H\semidirect U_{x}^{H}.
\end{equation*}
These local pictures can be thought of as two atlases of orbifold
charts on $|\widetilde{\grfont{G}}|$, endowing it with two smooth
orbifold structures. Indeed, one could use the local picture to
\emph{define} the topology on the fixed point sectors. However, in
nature one often encounters somewhat more general groupoids, and it is
worthwhile to have fixed point sectors defined intrinsically for them
as well. In fact, this will become crucial for us when we discuss
equivariant stable homotopy groups.

\subsection{Finer Structure}\label{subsec::finerstructure}

The fixed point sectors have a finer structure worth mentioning. Just
as the collection of fixed point sets of a $G$-space come equipped
with a web of maps making them into an
$\mathcal{O}(G)^{\text{op}}$-space, i.e., a contravariant diagram of
spaces over the \emph{orbit category} $\mathcal{O}(G)$, so also do our
sectors include into one another\footnote{One might also consider the
diagram over the \emph{subgroup category} of $G$, c.f.
\cite[p. 206]{luck;chern-characters-proper-equivariant-homology}.}. The
difference in our case is that one must restrict to components over a
given base point to ensure that the isotropy groups are actually
subconjugate in the groupoid. Also, the inclusions we obtain are in
general only defined in the localized category $\Orb$, rather than on
the level of groupoids. We will only discuss the situation for the
unreduced sectors $\widetilde{\grfont{G}}$; very similar results hold
for $\overline{\grfont{G}}$ as well.

By a \emph{connected component} or \emph{$\grfont{G}$-component} of a
groupoid $\grfont{G}$, we mean the inverse image of a component of
$|\grfont{G}|$ under the quotient map. Suppose that $\grfont{G}$ is a
proper \'etale groupoid. We wish to understand the connected
components of $\widetilde{\grfont{G}}$. Such a component corresponds
to a connected component of the classifying space, since arrows become
paths under realization. Consider the following equivalence relation
on $\widetilde{S}(\grfont{G})$. First, suppose $(p,H)$ and $(q,K)$ are
elements of $\widetilde{S}(\grfont{G})$ such that $p$ and $q$ lie in
the same equivariantly contractible chart $G_x\semidirect U_x$ about
some $x\in G_{0}$. Then (up to conjugacy) we may regard $H$ and $K$ as
subgroups of $G_{x}$. Write $(p,H)\equivalent_{\text{loc}} (q,K)$ if
and only if $(H)_{G_{x}} = (K)_{G_{x}}$. For general points in
$\widetilde{S}(\grfont{G})$, we let $(p,H)\equivalent (q,K)$ if there
is either a finite chain of $\equivalent_\text{loc}$'s joining them,
or else an arrow $g\in G_{1}$ with $g^{-1}Hg=K$. Equivalent elements
$(p,H)$ and $(q,K)$ will be called \emph{locally conjugate}.

The connected components of $\widetilde{\grfont{G}}$ are exactly these
local conjugacy classes. In the case of a global quotient
$G\semidirect M$, each sector corresponds to a disjoint union of fixed
point sets $G\semidirect \disunion_{K} M^{K}$ for $K\isom H$, and each
local conjugacy class corresponds to a $G$-component of this disjoint
union.

Let $p:\widetilde{\grfont{G}}\rightarrow \grfont{G}$ be the natural
projection. For each $x\in G_{0}$, the fiber $p^{-1}(x)\subseteq
\widetilde{\grfont{G}}$ is the discrete groupoid with objects the
subgroups of $G_{x}$, and morphisms given by the conjugation action of
$G_{x}$. For each subgroup $H\subseteq G_{x}$, let
\begin{equation*}
T_{H}\subseteq \widetilde{S}(\grfont{G})
\end{equation*}
be the connected component of $\widetilde{\grfont{G}}$ containing $(x,H)$;
so $T_{H}$ consists of the points locally conjugate to $(x,H)$. 

Now suppose $f:G_{x}/H\rightarrow G_{x}/K$ is an arrow in
$\mathcal{O}(G_{x})$. Then there is $g\in G_{x}$, determined up to
multiplication by an element of $K$ on the right, such that
$gHg^{-1}\subseteq K$ and $f(sH)=sgK$ for each $s\in G_{x}$. We define
an orbifold morphism $f^{\#}:\widetilde{\grfont{G}}|_{T_{K}}\rightarrow
\widetilde{\grfont{G}}|_{T_{H}}$. For $(x',K')\in T_{K}$, take a local
chart $G_{x'}\semidirect U_{x'}$ for $\grfont{G}$ at $x'$. Then
$\widetilde{\grfont{G}}|_{T_{K}}$ is locally isomorphic to
$G_{x'}\semidirect \disunion_{L\in (K')_{G_{x'}}} U_{x'}^{L}$ at $(x',
K')$. Similarly, $\widetilde{\grfont{G}}|_{T_{H}}$ is locally isomorphic
to $G_{x'}\semidirect \disunion_{N\in(H')_{G_{x'}}}U_{x'}^{N}$ at
$(x', H')$. Moreover, because $(x', K')$ is locally conjugate to $(x,
K)$ and $(x', H')$ is locally conjugate to $(x, H)$, we may identify
$H'$ with a subgroup of $K'$ up to conjugacy in $G_{x'}$. Fixing such
an identification determines a faithful homomorphism
\begin{equation*}
G_{x'}\semidirect \disunion_{L\in (K')_{G_{x'}}} U_{x'}^{L}\rightarrow
G_{x'}\semidirect \disunion_{N\in(H')_{G_{x'}}}U_{x'}^{N}.
\end{equation*}
Together, these local homomorphisms (defined up to conjugation by an
element of $G_{x'}$) determine the desired orbifold morphism
$f^{\#}$. As $f$ varies over $\mathcal{O}(G_{x})$, one obtains a
contravariant functor from $\mathcal{O}(G_{x})$ to $\Orb$. Note,
however, that all the morphisms $G_{x}/H\to G_{x}/K$ in
$\mathcal{O}(G_{x})$ go to the same orbifold morphism, since morphisms
that differ by conjugations in the local groups are identified in the
localized category $\Orb$. So the functor factors through the category
with objects the conjugacy classes of subgroups of $G_{x}$, and a
single morphism $(H)\rightarrow (K)$ if and only if $H$ is
subconjugate to $K$.

If $\grfont{G}$ is a general orbifold groupoid, we can arrange to have
a weak equivalence $\epsilon :\grfont{E}\rightarrow \grfont{G}$ where
$\grfont{E}$ is proper \'etale. In this case, the subgroupoids
generated by the images of the components $T_{H}$ give rise to a
Morita equivalent $\mathcal{O}(G_{x})^{\text{op}}$-orbifold involving
components of $\widetilde{\grfont{G}}$.

\section{Homotopy Groups}\label{sec::homotopyGroups}

Now that we know something about the fixed point data of an orbifold,
it is easy to define interesting new invariants. Here, we stick to
homotopy groups, but any other (weak) homotopy functor could also be
applied.

\subsection{Definitions}\label{subsec::definitionsHtpy}

The \emph{$n^{\text{th}}$ stable orbifold homotopy group} of an orbifold
groupoid $\grfont{G}$ is
\begin{equation}
\sthorb_{n}(\grfont{G}):=\stho_{n}(B\overline{\grfont{G}}_{+}),
\end{equation}
where on the right we take the ordinary stable homotopy of the
classifying space $B\overline{\grfont{G}}$ with a disjoint base point
$+$ added.

If $(x,H)$ is a base point in $\widetilde{\grfont{G}}$, the
\emph{$n^{\text{th}}$ (unstable) extended orbifold homotopy group} of
$\grfont{G}$ is
\begin{equation}
\piorbext_{n}(\grfont{G},(x,H)):=\pi_{n}(B\widetilde{\grfont{G}}, [x,H]),
\end{equation}
where $[x, H]$ is the point in the classifying space corresponding to
$(x,H)\in\widetilde{S}(\grfont{G})$.

When $\grfont{G}$ is a groupoid presentation of an orbifold
$\orbfont{X}$, we sometimes write $\sthorb (\orbfont{X})$ for
$\sthorb_{n}(\grfont{G})$ and $\piorb_{n}(\orbfont{X}, (x,H))$ for
$\piorb_{n}(\grfont{G}, (x,H))$. This is justified because both
constructions are Morita invariant: the weak homotopy types of
$B\widetilde{\grfont{G}}$ and $B\overline{\grfont{G}}$ depend only on
the Morita class of $\widetilde{\grfont{G}}$ and
$\overline{\grfont{G}}$, which in turn are determined by the Morita
class of $\grfont{G}$ by Lemma \ref{lem::sectors}. It follows that
both invariants descend to the localized category $\Orb$ of
orbifolds. We state this important fact as a theorem.

\begin{theorem}\label{thm::maintheorem}
The stable orbifold homotopy groups $\sthorb_{n}(\orbfont{X})$ and
extended orbifold homotopy groups $\piorb_{n}(\orbfont{X}, (x, H))$
are orbifold invariants. 
\end{theorem}

\begin{remark}
Regarding the ``finer structure'' discussed above, we see that
$\piorbext_{n}(\grfont{G},(x,-))$ may be viewed as an
$\mathcal{O}(G_{x})^{\text{op}}$-diagram of groups and (conjugacy
classes of) group homomorphisms as $H$ runs though the subgroups of
$G_{x}$.
\end{remark}

\subsection{Relation to Classical Theories}\label{subsec::classic}

To give some idea of what these new homotopy groups measure, we can
compare them with existing theories for orbifolds and equivariant
spaces. The main result is that our new stable homotopy groups include
equivariant stable homotopy groups as a special case.

\subsubsection{Classical Orbifold Homotopy Theory}
\label{subsubsec::classic-orbifold}

The extended orbifold homotopy groups are generalizations of the
classical orbifold homotopy groups. The latter appear as the
contribution of the untwisted sector $\widetilde{\grfont{G}}^{\langle
1\rangle}$ in the extended groups:
\begin{equation}
\piorb_{n}(\grfont{G}, x)\isom \piorbext_{n}(\grfont{G}, (x,\langle 1\rangle)).
\end{equation}
As we vary the group $H$ in the base point $(x, H)$, we have seen that
we obtain an $\mathcal{O}(G_{x})^{\text{op}}$-orbifold. For a given
$H$, the extended homotopy group corresponds to the classical orbifold
homotopy group $\piorb_{n}(\widetilde{\grfont{G}}^{H}, (x,H))$, which
will depend up to isomorphism only on the component $T_{H}$ of the
base point $(x,H)$. So in the equivariant situation
$\grfont{G}=G\semidirect M$, we are calculating the ordinary homotopy
groups of the Borel construction $ENH\times_{NH} M^{H}$. What's more,
Moerdijk has shown \cite{moerdijk;classifying-toposes-foliations} that
classes in $\piorb_{n}(\grfont{G}, x)$ can be represented by based
generalized morphisms $\bbS^{n}\to \grfont{G}$, where the sphere
$\bbS^{n}$ is viewed as a unit groupoid. Consequently, elements in
$\piorbext_{n}(\grfont{G}, (x,H))$ can be represented by pointed
generalized morphisms $\bbS^{n}\to\widetilde{\grfont{G}}^{H}$, which
themselves correspond to \emph{faithful} pointed generalized morphisms
$H\semidirect \bbS^{n}\to \grfont{G}$. Here, $H$ is acting trivially
on $\bbS^{n}$, so the translation groupoid is really the same thing as
the product $H\times \bbS^{n}$.

It is a remarkable fact in equivariant homotopy theory that a $G$-map
is an equivariant homotopy equivalence if and only if it induces
homotopy equivalences on all fixed point sets. It would be interesting
to see if any similar statement holds for orbifold groupoids.

\subsubsection{Equivariant Stable Homotopy Theory}
\label{subsubsec::equivar-stable}

It turns out that in good situations, the stable homotopy groups we
have defined for orbifolds really calculate a classical invariant. 

\begin{proposition}\label{prop::mainprop}
Let $\orbfont{X}=M/G$ where $M$ is a smooth manifold and $G$ is a
compact Lie group acting smoothly and almost freely. Then the total
stable equivariant homotopy group
\begin{equation*}
\stho^{G}_{\mathrm{tot}}(M_{+}):=\bigoplus_{n} \stho^{G}_{n}(M_{+})
\end{equation*}
is an orbifold invariant.
\end{proposition}
\begin{proof}
We need a calculational lemma.
\begin{lemma}\label{lem::calculation}
Let $\grfont{G}$ be an orbifold groupoid Morita equivalent to
$G\semidirect M$, where $G$ is a compact Lie group acting smoothly and
almost freely on the manifold $M$. Then 
\begin{equation*}
\sthorb_{n}(\grfont{G})\isom
\directsum_{(H)}\stho^{W_{G}H}_{n+d(h)}(EW_{G}H_{+}\smashprod M^{H}_{+}),
\end{equation*}
where $W_{G}H=N_{G}H/H$ and $d(H)=\dim_{\bbR}(W_{G}H)$.
\end{lemma}
\begin{proof}[Proof of Lemma:]
Since $\grfont{G}$ is Morita equivalent to $G\semidirect M$, their
fixed point sectors are also Morita equivalent by the Lemma
\ref{lem::sectors}. Also,
\begin{align*}
B\overline{\grfont{G}} &\weakhtpic 
B\left(\disunion_{(H)}W_{G}H\semidirect M^{H}\right)\\
          &\homotopic \disunion_{(H)} EW_{G}H\times_{W_{G}H} M^{H},
\end{align*}
where $\weakhtpic$ denotes weak homotopy equivalence. Thus, we
calculate:
\begin{align*}
\sthorb_{n}(\grfont{G})&\isom \stho_{n}((\disunion_{(H)}
EW_{G}H\times_{W_{G}H}M^{H})_{+})\\ 
&\isom \directsum_{(H)}\stho_{n}(EW_{G}H_{+}\smashprod_{W_{G}H}M^{H}_{+})\\
&\isom \directsum_{(H)}\stho_{n+d(H)}(\bbS^{d(H)}\smashprod EW_{G}H_{+}\smashprod_{W_{G}H}M^{H}_{+})\\
&\isom \directsum_{(H)}\stho_{n+d(H)}^{W_{G}H}(EW_{G}H_{+}\smashprod M^{H}_{+}).
\end{align*}
The last isomorphism is a Wirthm\"uller isomorphism, obtained by
noting that $EW_{G}H_{+}\smashprod M^{H}_{+}$ is a free $W_{G}H$-space. Here,
we have regarded $\bbS^{d(H)}$ as the representation sphere corresponding to
the adjoint action of the trivial subgroup of $W_{G}H$ on $T_{e}W_{G}H$. 
\end{proof}
The proposition now follows, since the lemma shows that the total stable
equivariant homotopy of $M_{+}$ is isomorphic to the total orbifold
stable homotopy of $\grfont{G}$, which is an orbifold invariant by
Theorem \ref{thm::maintheorem}. 
\end{proof}

\begin{remark}
The above is the strongest result that could be hoped for while
maintaining Morita invariance. The numbers $d(H)$ are manifestly not
Morita invariant, for given $G\semidirect M$, the groupoid $(G\times
\bbS^{1})\semidirect (M\times \bbS^{1})$ is Morita equivalent but has all
$d(H)$ increased by one.
\end{remark}

\begin{example}
Let $\orbfont{X}$ be an \emph{effective} $n$-dimensional
orbifold. Then $\orbfont{X}$ may be presented as a quotient $F/O(n)$,
where $F$, the \emph{frame bundle} of $\orbfont{X}$, is a smooth
manifold with an almost free $O(n)$-action. Then 
\begin{equation*}
\stho_{\mathrm{tot}}^{O(n)}(F_{+})\isom \sthorb_{\mathrm{tot}}(\orbfont{X})
\end{equation*}
is an orbifold invariant.
\end{example}

We can draw another interesting conclusion from the lemma.

\begin{corollary}\label{cor::maincor}
If $M/G$ and $M'/G'$ are two global quotient presentations of the same
orbifold $\orbfont{X}$, then there are isomorphisms
\begin{equation*}
\xymatrix{ 
\stho^{G}_{n}(M_{+}) & \sthorb_{n}(\orbfont{X})
\ar[l]_-{\isom}\ar[r]^-{\isom} & \stho_{n}^{G'}(M_{+}') 
}
\end{equation*}
for each integer $n$.
\end{corollary}
\begin{proof}
In this case, $d(H)$ is zero for every subgroup of $G$ or $G'$. The
result follows immediately from tom Dieck's isomorphism (Equation
\eqref{eqn::tomDieck}).
\end{proof}

\begin{example}
Suppose $\orbfont{X}$ is a connected global quotient orbifold. Then
the orbifold universal cover (c.f. \cite{thurston-notes}) of
$\orbfont{X}$ is a manifold $Y$. The orbifold fundamental group
$\Gamma:=\piorb_{1}(\orbfont{X})$ acts on $Y$ with quotient
$\orbfont{X}$. If $H_{1}$ and $H_{2}$ are any two normal subgroups of
$\Gamma$ acting freely on $Y$, then $M_{1}/G_{1}$ and $M_{2}/G_{2}$
are two different global quotient presentations of $\orbfont{X}$,
where $M_{i}=Y/H_{i}$ and $G_{i}=\Gamma/H_{i}$. In fact, any two such
presentations where the $M_{i}$ are connected arise in this way. The
corollary says that $\stho_{n}^{G_{1}}(M_{1+})\isom
\stho_{n}^{G_{2}}(M_{2+})$, which may be confirmed via explicit
calculation.
\end{example}

\section{Further Questions}\label{sec::future}

We have called our new invariants ``homotopy groups,'' but it is still
unclear exactly what is meant by homotopy equivalent orbifold
groupoids or orbifolds. Natural transformations between groupoid
homomorphisms realize to homotopies of maps between classifying
spaces; so our groups have at least one sort of homotopy
invariance. However, the interval $I$ can be viewed as a unit
groupoid, and so one could also study homotopies of the form
$\grfont{G}\times I\rightarrow \grfont{H}$. These realize into maps
$B\grfont{G}\times (\Delta^{\infty}\times I)\rightarrow B\grfont{H}$,
whose significance is less clear.

A more comprehensive treatment would involve setting up a model
structure on the category $\Gpd$ (or $\Orb$) and studying its
relationship to our homotopy group functors. In a future paper, we
hope to construct such a model category using techniques related to
the homotopy theory of schemes and the universal homotopy theory of
Dugger \cite{dugger;universal-homotopy-theories}, as suggested by Dan
Isaksen \cite{isaksen;private-communication}. The motivation for this
approach is as follows: our fixed point sectors can be identified with
mapping spaces of (faithful) homomorphisms from finite groups into
$\grfont{G}$. In broad strokes, looking at morphisms from various
groups into the groupoid $\grfont{G}$ is analogous to studying various
sorts of points in a scheme or stack. Indeed, this should already be a
familiar notion for equivariant topologists, given that one often
identifies $M^{H}$ with the mapping space $\map_{G}(G/H,M)$. The hope
is that this abstract approach might shed light on the question at the
end of Section \ref{subsubsec::classic-orbifold}: namely, to what
extent the classical orbifold homotopy types of the fixed point
sectors determine the ``homotopy type'' of the orbifold in some
reasonable model structure?  Progress along similar lines appears in
the very recent preprint of Noohi
\cite{noohi;homotopy-topological-stacks-I} regarding homotopy groups
of topological stacks (see also \cite{jardine;simplicial-presheaves},
\cite{hollander;homotopy-theory-stacks}).

\appendix 

\section{Classifying Spaces of Translation Groupoids}
\label{app::borelconstruction}

Let $G$ be a compact Lie group acting on $M$, and let
$\grfont{G}=G\semidirect M$ be the translation groupoid. The homotopy
equivalence $B\grfont{G}\homotopic EG\times_{G}M$ seems to be a folk
theorem. We give a quick proof generalizing Segal's arguments for the
group case in \cite{segal;classifying-spaces-spectral-sequences}.

Let $\Gbar$ be the groupoid with objects $G$ and arrows $G\times G$,
so that there is a unique isomorphism $(g_{1}, g_{2}):g_{1}\to g_{2}$
between every two objects. Hence, this category is equivalent to the
trivial category with one object and one morphism. So $B\Gbar$ is a
model for $EG$, since the free simplicial $G$-action on the nerve
realizes to a free $G$-action on the contractible space $B\Gbar$. Now
consider the product groupoid $\Gbar\times M$, where $M$ denotes the
unit groupoid on $M$. Since $B$ respects products, $B(\Gbar\times
M)=B\Gbar\times BM\homotopic EG\times M$.

$G$ acts on the groupoid $\Gbar\times M$ by automorphisms as
follows. For each object $(g_{1}, m)$, let $g(g_{1}, m)=(g_{1}g^{-1},
gm)$; for each arrow $(g_{1}, g_{2}, m)$, let $g(g_{1}, g_{2},
m)=(g_{1}g^{-1}, g_{2}g^{-1}, gm)$. This action induces a simplicial
$G$-action on the nerve, and one obtains the diagonal action on
$EG\times M$ upon realization. Moreover, the quotient groupoid under
the action is isomorphic to $G\semidirect M$ via the homomorphism
sending the object orbit $G(g_{1}, m)$ to $g_{1}m$ and the arrow orbit
$G(g_{1}, g_{2}, m)$ to $(g_{2}g_{1}^{-1}, g_{1}m)$. Consequently, we
may identify $EG\times_{G} M \homotopic B(\Gbar\times M)/G$ with
$B(G\semidirect M)$, up to homotopy.

\bibliography{orbifolds}
\bibliographystyle{plain}

\end{document}